\documentclass[12pt]{amsart}
\usepackage{amsmath,amssymb}
\usepackage{amsfonts}
\usepackage{amsthm}
\usepackage{latexsym}
\usepackage{graphicx}

\def\lf{\left}
\def\ri{\right}

\def\p{\partial}

\def\R{\mathbb{R}}

\def\vv<#1>{\langle#1\rangle}

\def\bn{\mbox{\bf{n}}}

\def\XXint#1#2{\setbox0=\hbox{$#1{#2}{\int}$}{#2}\kern-.5\wd0 }

\def\XXint#1#2#3{{\setbox0=\hbox{$#1{#2#3}{\int}$}
     \vcenter{\hbox{$#2#3$}}\kern-.5\wd0}}

\def\vv<#1>{{\left\langle#1\right\rangle}}

\def\bt{\mbox{\bf{t}}}
\newtheorem{thm}{Theorem}[section]

\theoremstyle{definition}

\theoremstyle{remark}

\numberwithin{equation}{section}

\begin{document}
\title{A higher dimensional generalization of Hersch-Payne-Schiffer inequality for Steklov eigenvalues }

\author{Liangwei Yang}
\address{Department of Mathematics, Shantou University, Shantou, Guangdong, 515063, China}
\email{13lwyang@stu.edu.cn}
\author{Chengjie Yu$^1$}
\address{Department of Mathematics, Shantou University, Shantou, Guangdong, 515063, China}
\email{cjyu@stu.edu.cn}
\thanks{$^1$Research partially supported by a supporting project from the Department of Education of Guangdong Province with contract no. Yq2013073, and NSFC 11571215.}
\renewcommand{\subjclassname}{%
  \textup{2010} Mathematics Subject Classification}
\date{}
\begin{abstract}
In this paper, we generalize the Hersch-Payne-Schiffer inequality for Steklov eigenvalues to higher dimensional case by extending the trick used by Hersch, Payne and Schiffer to higher dimensional manifolds.
\end{abstract}
\maketitle\markboth{Liang \& Yu}{Hersch-Payne-Schiffer Inequality}
\section{Introduction}
Let $(M^n,g)$ be a compact connected Riemannian manifold with nonempty boundary. The Dirichlet-to-Neumann map $L:C^\infty(\partial M)\to C^\infty(\partial M)$ for functions is defined as
$$L(u)=\frac{\partial \hat u}{\partial \nu},$$
where $\hat u$ is the harmonic extension of $u\in C^\infty(\partial M)$ and $\nu$ is the outward normal vector field on $\partial M$. It is shown in \cite{Ta} that $L$ is a first order nonnegative self-adjoint elliptic pseudo-differential operator. So the spectrum of $L$ is discrete and can be arranged in increasing order (counting multiplicities) as:
$$0=\sigma_0<\sigma_1\leq \sigma_2\leq \cdots.$$
$\sigma_k$ is called the $k$-th Steklov eigenvalue of $(M^n,g)$. The Steklov eigenvalues have been intensively studied recently. \cite{GP} makes an excellent survey for recent progresses of the topic.

In 1974, Hersch,Payne and Schiffer \cite{HPS} proved the following inequality:
\begin{equation}\label{eqn-HPS}
\sigma_p(\Omega)\sigma_q(\Omega)L(\partial \Omega)^2\leq\left\{\begin{array}{ll}(p+q-1)^2\pi^2&\mbox{if $p+q$ is odd}\\
(p+q)^2\pi^2&\mbox{if $p+q$ is even}
\end{array}\right.
\end{equation}
by an elegant trick using conjugate harmonic functions. Here $\Omega$ is a bounded simply connected domain in $\R^2$ with smooth boundary and $L(\partial \Omega)$ is the length of $\partial \Omega$. In fact, in \cite{HPS}, Hersch, Payne and Schiffer also obtained similar inequalities for multiple connected domains.

The Hersch-Payne-Schiffer inequalities was recently generalized to surfaces of higher genus by  Girouard and Polterovich  via different techniques using the Ahlfors map in \cite{GP1}. In fact, they proved:
\begin{equation}\label{eqn-GP}
\sigma_p(\Omega)\sigma_q(\Omega)L(\partial \Omega)^2\leq\left\{\begin{array}{ll}(p+q-1)^2(\gamma+l)^2\pi^2&\mbox{if $p+q$ is odd}\\
(p+q)^2(\gamma+l)^2\pi^2&\mbox{if $p+q$ is even.}
\end{array}\right.
\end{equation}
Here $\Omega$ is a compact oriented surface with genus $\gamma$ and $l$ connected boundary components. When $\gamma=0$ and $l=1$, it recovers the Hersch-Payne-Schiffer inequality for bounded simply connected plane domains. However, when $\gamma=0$ and $l>1$, \eqref{eqn-GP} is different with the Hersch-Payne-Schiffer inequality for multiple connected plane domains in \cite{HPS}.

In \cite{RS2}, Raulot and Savo introduced the following Dirichlet-to-Neumann map for differential forms. Let $\omega\in A^r(\partial M)$, the Dirichlet-to-Neumann map $L^{(r)}$ for $r$-differential forms is defined as
$$L^{(r)}(\omega)=i_{\nu}d\hat\omega.$$
Here $\hat\omega$ is the tangential harmonic extension of $\omega$. It is shown in \cite{RS2} that $L^{(r)}$ is also a first order nonnegative self-adjoint elliptic pseudo-differential operator. So, the spectrum of $L^{(r)}$ is discrete. We arrange the eigenvalues of $L^{(r)}$ in increasing order (counting multiplicities) as:
\begin{equation}
0\leq \sigma_1^{(r)}\leq \sigma_2^{(r)}\leq\cdots.
\end{equation}
There is a different Dirichlet-to-Neumann map for differential forms in \cite{BS} where the authors used another kind of harmonic extension for differential forms on the boundary. However, the Dirichlet-to-Neumann map introduced in \cite{BS} is not elliptic for differential forms with positive degree, because the kernel of the map is not finite dimensional.

In this paper, by extending the trick of Hersch-Payne-Schiffer to higher dimensional cases, we obtain the following generalization of the Hersch-Payne-Schiffer inequality.
\begin{thm}\label{thm-main}
Let $(M^n,g)$ be a compact oriented Riemannian manifold with nonempty boundary. Let $\sigma_k^{(r)}(M)$ be the $k$-th Steklov eigenvalue for differential $r$-forms of $\partial M$ and $\lambda_k(\p M)$ be the $k$-th eigenvalue for the Hodge-Laplacian operator of $\p M$ (both counting multiplicities).
Then, for any two positive integers $p$ and $q$, we have
\begin{equation}
\sigma_p^{(0)}(M)\sigma_{b_{n-2}+q}^{(n-2)}(M)\leq \lambda_{p+q+b_{n-1}-1}(\p M)
\end{equation}
where $b_{k}$ is the $k$-th Betti number of $M$.
\end{thm}

When $M$ is a compact oriented surface with genus $0$ and $l$ connected boundary components, Theorem \ref{thm-main} is just the Hersch-Payne-Schiffer inequality \eqref{eqn-HPS}. When $M$ is a compact oriented surface with genus $\gamma$ and $l$ connected boundary components, Theorem \ref{thm-main} is a higher genus generalization of the Hersch-Payne-Schiffer inequality different with \eqref{eqn-GP}. For example, when $M$ is a compact oriented surface of genus $\gamma$ with connected boundary, Theorem \ref{thm-main} give us the following inequality
directly:
\begin{equation}
\sigma_p(M)\sigma_q(M)L(\partial M)^2\leq \left\{\begin{array}{ll}(p+q+2\gamma-1)^2\pi^2&\mbox{if $p+q$ is odd}\\
(p+q+2\gamma)^2\pi^2&\mbox{if $p+q$ is even.}
\end{array}\right.
\end{equation}

In \cite{K,RS1,RS2,RS3}, the authors also obtained interesting estimates of  the Steklov eigenvalues for differential forms under some nonnegative assumptions on the curvature of the manifold and its boundary. Comparing to these results, our result is free of curvature assumptions. It is just a direct generalization of the Hersch-Payne-Schiffer inequality in \cite{HPS}.
\section{Preliminaries}
In this section, we recall some preliminaries in Hodge Theory for oriented compact Riemannian manifolds with nonempty boundary.

Let $(M^n,g)$ be an $n$-dimensional oriented compact Riemannian manifold with nonempty boundary and $*$ be the Hodge star operator and
$$\delta=(-1)^{nr+1}*d*:A^{r+1}(M)\to A^r(M)$$
be the formal adjoint of $d:A^r(M)\to A^{r+1}(M)$. Here $A^r(M)$ is the space of differential $r$-forms on $M$. Let $i:\p M\to M$ be natural embedding, $\bt \omega=i^*\omega$ for any differential forms $\omega$ on $M$, and
$$\bn \omega=\omega-\bt\omega$$
on $\p M$. The following relations of $*$, $\bt$ and $\bn$ are useful:
\begin{equation}
\bt*\omega=*\bn \omega,\ \ \bn*\omega=*\bt\omega
\end{equation}
for any differential form $\omega$. By Stokes formula and the relation above, we have
\begin{equation}\label{eqn-Green}
\begin{split}
\int_M\vv<d\alpha,\beta>dV=&\int_M\vv<\alpha,\delta\beta>dV+\int_{\p M}\bt\alpha\wedge*\bn\beta\\
=&\int_M\vv<\alpha,\delta\beta>dV+\int_{\p M}\vv<i^*\alpha,i_\nu\beta>dA.
\end{split}
\end{equation}
Here $\nu$ is the unit outward normal of $M$ and
$$i_\nu\beta(\cdot)=\beta(\nu,\cdot).$$

The following Hodge-Freiderich-Morrey decomposition can be found in \cite{S}.
\begin{thm}\label{thm-Hodge}
Let $(M^n,g)$ be an $n$-dimensional oriented compact Riemannian manifold with nonempty boundary. Let
$$\mathcal{E}^r(M)=\{d\alpha\ |\ \alpha\in A^{r-1}(M), \bt \alpha=0\},$$
$$\mathcal{C}^r(M)=\{\delta\beta\ |\ \beta\in A^{r+1}(M), \bn \beta=0\}$$
and
$$\mathcal{H}^r(M)=\{\gamma\in A^r(M)\ |\ d\gamma=0,\ \mbox{and}\ \delta\gamma=0\}.$$
Then, we have the following $L^2$-orthogonal decomposition of $A^r(M)$:
$$A^r(M)=\mathcal{E}^r(M)\oplus\mathcal{C}^r(M)\oplus\mathcal{H}^r(M),$$
for any $r=0,1,\cdots,n$. Moreover, we have the following $L^2$-orthogonal decompositions for $\mathcal H^r(M)$:
$$\mathcal H^r(M)=\mathcal H_D^r(M)\oplus \mathcal H^r_{co}(M)$$
and
$$\mathcal H^r(M)=\mathcal H_N^r(M)\oplus \mathcal H^r_{ex}(M)$$
where
$$H_D^r(M)=\{\gamma\in \mathcal H^r(M)\ |\ \bt\gamma=0\},$$
$$H_{co}^r(M)=\{\gamma\in \mathcal H^r(M)\ |\ \gamma=\delta\alpha\ \mbox{for some}\ \alpha\in A^{r+1}(M)\},$$
$$\mathcal H^r_N(M)=\{\gamma\in \mathcal H^r(M)\ |\ \bn\gamma=0\}$$
and
$$\mathcal H^r_{ex}(M)=\{\gamma\in \mathcal H^r(M)\ |\ \gamma=d\alpha\ \mbox{for some}\ \alpha\in A^{r-1}(M)\}.$$
\end{thm}
Let $\alpha\in A^{r}(\p M)$ and $\hat \alpha\in A^r(M)$ be the unique tangential harmonic extension of $\alpha$. That is,
\begin{equation}
\left\{\begin{array}{l}\Delta \hat\alpha=0\\
\bt \hat \alpha=\alpha\\
\bn\hat\alpha=0
\end{array}\right.
\end{equation}
Here $\Delta=d\delta+\delta d$ is the Hodge-Laplacian operator. The
Dirichlet-to-Neumann map $L^{(r)}:A^r(\p M)\to A^r(\p M)$ for differential $r$-forms is defined as
\begin{equation}
L^{(r)}(\alpha)=i_\nu d\hat\alpha.
\end{equation}

It is proved in \cite{RS2} that $L^{(r)}$ is a nonnegative elliptic self-adjoint pseudo-differential operator of first order. So, the eigenvalues of $L^{(r)}$ are discrete which are called the Steklov eigenvalues for differential $r$-forms. By Theorem \ref{thm-Hodge} and \eqref{eqn-Green}, it is clear that
$$\dim \ker L^{(r)}=\dim \mathcal H_N^r(M)=b_r$$
where $b_r$ is the $r$-th Betti number of $M$. So, the multiplicity of the eigenvalue $0$ is $b_r$ for $L^{(r)}$. Arrange all the eigenvalues of $L^{(r)}$ (counting multiplicities) as follows:
$$0\leq\sigma_1^{(r)}\leq\sigma_2^{(r)}\leq\cdots.$$
It is clear that $\sigma^{(0)}_1,\sigma^{(0)}_2,\cdots$ are the Steklov eigenvalues in usual sense. Let $\alpha_k^{(r)}\in A^{r}(\p M)$ be the eigenform for $\sigma^{(r)}_k$. Then, we have the following variational principle for $\sigma^{(r)}_k$:
\begin{equation}
\sigma^{(r)}_k=\inf\lf\{\frac{\int_M\vv<d\hat\alpha,d\hat\alpha>+\vv<\delta\hat \alpha,\delta\hat\alpha>dV}{\int_{\p M}\vv<\alpha,\alpha>dA}\ \Bigg|\begin{array}{l}\alpha\in A^{k}(\partial M),\alpha\neq 0,\\
\alpha\perp \alpha_1^{(r)},\cdots,\alpha\perp \alpha_{k-1}^{(r)}\end{array} \ri\}.
\end{equation}
\section{Proof of Theorem \ref{thm-main}}
We are now ready to prove Theorem \ref{thm-main}.
\begin{proof}
Let $\phi_k\in C^\infty(\p M)$ be the eigenfunction of $\lambda_k(\p M)$. Since $\lambda_1(\partial M)=0$. We can suppose that $\phi_1\equiv 1$ on $\partial M$. Let
\begin{equation}\label{eqn-u}
u=c_2\phi_2+c_2\phi_3+\cdots+c_{p+q+b_{n-1}-1}\phi_{p+q+b_{n-1}-1}\in C^{\infty}(\p M)
\end{equation}
where $c_2,c_3,\cdots,c_{p+q+b_{n-1}-1}$ are constants that are not all zero to be determined. Then, $u$ is not a constant function on $\partial M$. Note that
\begin{equation}
d*d\hat u=*\delta d\hat u=*\Delta\hat u=0.
\end{equation}
So, $*d\hat u$ is a closed $(n-1)$-form on $M$. If we choose the constants in \eqref{eqn-u} such that
\begin{equation}\label{eqn-require-1}
*d\hat u\perp \{\gamma\in \mathcal H^{n-1}(M)\ |\ \bn \gamma=0\}
\end{equation}
Then, by Theorem \ref{thm-Hodge}, we know that $*d\hat u$ is exact. Suppose
\begin{equation}
*d\hat u=d\tilde \alpha
\end{equation}
with $\tilde \alpha\in A^{n-2}(M)$. By Theorem \ref{thm-Hodge} again, suppose that
\begin{equation}
\tilde\alpha=d\eta+\delta\beta+\gamma
\end{equation}
where $\eta\in A^{n-3}(M)$, $\beta\in A^{n-1}(M)$ with $\bn \beta=0$ and $\gamma\in \mathcal H^{n-2}(M)$ with $\bn\gamma=0$. Note that
\begin{equation}
\ker L^{(n-2)}=\{i^*\omega\ |\ \omega\in \mathcal H^{n-2}(M),\ \bn\omega=0\}.
\end{equation}
So, we can choose a suitable $\tilde \gamma\in \mathcal H^{n-2}(M)$ with $\bn\tilde \gamma=0$ such that
\begin{equation}
i^*(\delta\beta+\tilde \gamma)\perp \ker L^{(n-2)}.
\end{equation}
Let $\alpha=\delta\beta+\tilde \gamma$. Then
\begin{equation}
d\alpha=d\delta\beta=d\tilde\alpha=*d\hat u
\end{equation}
and
\begin{equation}
\delta\alpha=0.
\end{equation}
Moreover,
\begin{equation}
\Delta\alpha=(d\delta+\delta d)\alpha=\delta d\alpha=\delta*d\hat u=(-1)^{n-1}*dd\hat u=0,
\end{equation}
and
\begin{equation}
\bn \alpha=\bn\delta\beta+\bn\tilde\gamma=\delta\bn\beta=0.
\end{equation}
So, $\alpha$ can be view as the tangential harmonic extension of $i^*\alpha\in A^{n-2}(\p M).$ Furthermore, we choose the constants in \eqref{eqn-u}, so that
\begin{equation}\label{eqn-require-2}
  i^*\alpha\perp \alpha_{b_{n-2}+1}^{(n-2)},\cdots,\alpha_{b_{n-2}+q-1}^{(n-2)},
\end{equation}
and
\begin{equation}\label{eqn-require-3}
u\perp \alpha_2^{(0)},\cdots,\alpha_{p-1}^{(0)}.
\end{equation}
The requirements \eqref{eqn-require-1},\eqref{eqn-require-2} and \eqref{eqn-require-3} form a homogeneous linear system of $c_2,c_3,\cdots,c_{p+q+b_{n-1}-1}$ with $p+q+b_{n-1}-3$ equations. So, there are constants $c_2,c_3,\cdots,c_{p+q+b_{n-1}-1}$ that are not all zero satisfying the requirements \eqref{eqn-require-1},\eqref{eqn-require-2} and \eqref{eqn-require-3}.

Note that, by \eqref{eqn-Green},
\begin{equation}\label{eqn-da}
\begin{split}
&\int_M\vv<d\alpha,d\alpha>dV\\
=&\int_{\p M}\vv<i^*\alpha,i_\nu d\alpha>dA\\
\leq &\lf(\int_{\p M}\vv<i^*\alpha,i^*\alpha>dA\ri)^{1/2}\lf(\int_{\p M}\vv<i_\nu d\alpha,i_\nu d\alpha>dA\ri)^{1/2}\\
=&\lf(\int_{\p M}\vv<i^*\alpha,i^*\alpha>dA\ri)^{1/2}\lf(\int_{\p M}\vv<i_\nu *d\hat u,i_\nu *d\hat u>dA\ri)^{1/2}\\
=&\lf(\int_{\p M}\vv<i^*\alpha,i^*\alpha>dA\ri)^{1/2}\lf(\int_{\p M}\vv<du,d u>dA\ri)^{1/2}.\\
\end{split}
\end{equation}
Moreover,
\begin{equation}\label{eqn-du}
\int_M\vv<d\hat u,d\hat u>dV=\int_M\vv<*d\hat u,*d\hat u>dV=\int_M\vv<d\alpha,d\alpha>dV.
\end{equation}
Therefore, by \eqref{eqn-da}, \eqref{eqn-du} and \eqref{eqn-u},
\begin{equation}
\begin{split}
&\sigma_p^{(0)}(M)\sigma_{b_{n-2}+q}^{(n-2)}(M)\\
\leq&\frac{\int_M\vv<d\hat u,d\hat u>dV}{\int_{\p M}u^2dA}\frac{\int_M\vv<d\alpha,d\alpha>dV}{\int_{\p M}\vv<i^*\alpha,i^*\alpha>dA}\\
=&\frac{\int_M\vv<d\alpha,d\alpha>dV}{\int_{\p M}u^2dA}\frac{\int_M\vv<d\alpha,d\alpha>dV}{\int_{\p M}\vv<i^*\alpha,i^*\alpha>dA}\\
\leq&\frac{\int_{\p M}\vv<du,du>dA}{\int_{\p M}u^2dA}\\
\leq&\lambda_{p+q+b_{n-1}-1}(\p M).
\end{split}
\end{equation}
\end{proof}


\begin{thebibliography}{99}
\bibitem{BS}Belishev, M.; Sharafutdinov, V. {\it Dirichlet to Neumann operator on differential forms.} Bull. Sci. Math. 132 (2008), no. 2, 128--145.
\bibitem{GP}Girouard, A.; Polterovich, I. {\it Spectral geometry of the Steklov problem.}  arXiv:1411.6567.
\bibitem{GP1}Girouard, A.; Polterovich, I. {\it Upper bounds for Steklov eigenvalues on surfaces.} Electron. Res. Announc. Math. Sci. 19 (2012), 77--85.
\bibitem{HPS} Hersch, J.; Payne, L. E.; Schiffer, M. M. {\it Some inequalities for Stekloff eigenvalues.} Arch. Rational Mech. Anal. 57 (1975), 99--114.
\bibitem{K}Kwong, K.-K. {\it Some sharp eigenvalue estimate for differential forms.} Private communication.
\bibitem{RS1}Raulot, S.; Savo, A. {\it On the spectrum of the Dirichlet-to-Neumann operator acting on forms of a Euclidean domain.} J. Geom. Phys. 77 (2014), 1--12.
\bibitem{RS2}Raulot, S.; Savo, A. {\it On the first eigenvalue of the Dirichlet-to-Neumann operator on forms.} J. Funct. Anal. 262 (2012), no. 3, 889--914.
\bibitem{RS3}Raulot, S.; Savo, A. {\it A Reilly formula and eigenvalue estimates for differential forms.} J. Geom. Anal. 21 (2011), no. 3, 620¨C640.
\bibitem{S}Schwarz, G. {\it Hodge decomposition¡ªa method for solving boundary value problems.} Lecture Notes in Mathematics, 1607. Springer-Verlag, Berlin, 1995. viii+155 pp. ISBN: 3-540-60016-7
\bibitem{Ta}Taylor, Michael E. {\it Partial differential equations II. Qualitative studies of linear equations.} Second edition. Applied Mathematical Sciences, 116. Springer, New York, 2011. xxii+614 pp. ISBN: 978-1-4419-7051-0.

\end{thebibliography}
\end{document}